\documentclass[12pt]{article}
\usepackage{t1enc}
\usepackage[latin1]{inputenc}
\usepackage[english]{babel}
\usepackage{latexsym}
\newcommand{\bp}{{\bf P}}
\newtheorem{propo}{Proposition}[section]

\newcommand{\be}{{\bf E}}
\newcommand{\bee}{{\bf e}}
\newtheorem{theorem}{Theorem}[section]

\newtheorem{corollary}{Corollary}[section]
\newtheorem{lemma}{Lemma}[section]

\def\kz{{\bf z}}
\def\bo{{\bf 0}}
\def\zd{{\cal Z}_d}
\def\z2{{\cal Z}_2}
\def\begg{\begin{equation}}
\def\endd{\end{equation}}
\def\bege{\begin{eqnarray}}
\def\ende{\end{eqnarray}}
\def\pe{{\bf P}}
\def\ep{{\varepsilon}}

\begin{document}

\centerline{\Large\bf On the local time of the asymmetric
Bernoulli walk}
\medskip
\bigskip\bigskip
\centerline{\it Dedicated to Professor S\'andor Cs\"org\H o on his
sixtieth birthday}

\bigskip \bigskip \bigskip \bigskip \bigskip

\renewcommand{\thefootnote}{1} \noindent
\textbf{Endre Cs\'{a}ki}\footnote{Research supported by the Hungarian
National Foundation for Scientif\/ic Research, Grant No.
K 61052 and K 67961.}\newline
Alfr\'ed R\'enyi Institute of Mathematics, Hungarian Academy of Sciences,
Budapest, P.O.B. 127, H-1364, Hungary. E-mail address: csaki@renyi.hu

\bigskip

\renewcommand{\thefootnote}{2} \noindent \textbf{Ant\'{o}nia
F\"{o}ldes}\footnote{Research supported by a PSC CUNY Grant,
No. 66494-0035.}\newline
Department of Mathematics, College of Staten Island, CUNY,
2800 Victory Blvd., Staten Island, New York
10314, U.S.A. E-mail address: foldes@mail.csi.cuny.edu

\bigskip

\noindent \textbf{P\'al R\'ev\'esz}$^1$ \newline
Institut f\"ur Statistik und Wahrscheinlichkeitstheorie, Technische
Universit\"at Wien, Wiedner Hauptstrasse 8-10/107 A-1040 Vienna, Austria.
E-mail address: reveszp@renyi.hu

\bigskip \bigskip \bigskip

\noindent \textit{Abstract}: We study some properties of the
local time of the asymmetric Bernoulli walk on the line. These
properties are very similar to the corresponding ones of the
simple symmetric random walks in higher ($d\geq3$) dimension,
which we established in  the recent years. The goal of this paper
is to highlight these similarities.

\bigskip


\bigskip

 \vspace{.1cm}

\vfill
\renewcommand{\thesection}{\arabic{section}.}
\section{Introduction}

\renewcommand{\thesection}{\arabic{section}} \setcounter{equation}{0}
\setcounter{theorem}{0} \setcounter{lemma}{0} 

The study of the local time (number of visits) of transient random walks 
started with the landmark papers of  Dvoretzky and Erd\H os \cite{DE50}
and Erd\H os and Taylor \cite{ET60}, who investigated the
properties of simple symmetric random walk in dimension $d\geq 3$,
in which case the random walk is transient by P\'olya theorem. For
further results we mention the books by Feller \cite{F1},
\cite{F2}, Spitzer \cite{S76} and R\'ev\'esz \cite{R05}.  In the recent
years our investigations were concentrated  on some of the fine
properties of the local and occupation times  of these walks. It is
well known that the simple asymmetric Bernoulli walk on the line is also
transient and as such it behaves similarly to other transient walks.
The goal of the present paper to put into evidence that many of the fine
properties of the local and occupation times  which we studied
for the $d(\geq3)$ dimensional transient symmetric walks are
really shared by  the asymmetric one dimensional Bernoulli walk.

Here we would like to discuss  only  the following three major topics.
\begin{itemize}
\item{} limit theorems for multiple visited points
\item{} joint behavior of local and occupation times
\item{} the local time around frequently visited points
\end{itemize}
These results in higher dimension were presented in our papers
\cite{CsFR06}, \cite{CsFR1}, \cite{CsFR2}.

In \cite{CsFR3}, a recent survey paper on these topics some of our
present results were given without proof. In this paper we would like
to collect all the results which we are having so far on the
asymmetric Bernoulli walk. Clearly to give full proofs for all
these results are very tedious but saying only, that proofs are
similar to the symmetric $d$-dimensional walk case is unfair. So we
take the middle way, namely we give   some of the proofs with an
emphasis on the differences between the two situations.

The organization of the paper is as follows. In Section 2 we collect
the relevant results on the simple symmetric $d$-dimensional walk. In
Section 3 we present the new results on  the local and occupation times
of the asymmetric one dimensional Bernoulli walk. In Section 4 we will
present some lemmas needed later in the proofs. In Section 5 we prove
Theorem 3.2. The proofs of Theorems 3.3 and 3.4 will be given in Section
6, while in Section 7 the proofs of Theorems 3.5 and 3.6 are presented.
Finally, Section 8 contains some remarks.

\renewcommand{\thesection}{\arabic{section}.}

\section{Random walk in higher dimension}

\renewcommand{\thesection}{\arabic{section}} \setcounter{equation}{0}
\setcounter{theorem}{0} \setcounter{lemma}{0}  
Let $\{\mathbf{S}_n\}_{n=1}^{\infty}$ be a symmetric random walk
starting at the origin $\mathbf{0}$ on the $d$%
-dimensional integer lattice $\mathcal{Z}_d$ where  $d\geq3$,  i.e.
$\mathbf{S}_0=\mathbf{0}$,
$\mathbf{S}_n=\sum_{k=1}^{n} \mathbf{X}_k$, $n=1,2,\dots$, where $\mathbf{X%
}_k,\, k=1,2,\dots$ are i.i.d. random variables with distribution
\[
\mathbf{P} (\mathbf{X}_1=\mathbf{e}_i)=\frac{1}{2d},\qquad
i=1,2,\ldots,2d
\]
and $\{\mathbf{e}_1,\mathbf{e}_2,...\mathbf{e}_d\}$ is a system of
orthogonal unit vectors in $\mathcal{Z}_d$ and $\mathbf{e}_{d+j}=-\mathbf{e}%
_j,$ $j=1,2,\ldots,d.$ Define the local time of the walk by
\begin{equation}
\xi(\mathbf{z},n):=\#\{k: \,\,0< k \leq n,\,\,\, \mathbf{S}
_k=\mathbf{z} \},\quad n=1,2,\ldots,  \label{loc1}
\end{equation}
where $\mathbf{z}$ is any lattice point of $\mathcal{Z}_d.$ 
Let $\xi(\mathbf{z},\infty):=\lim_{n\to\infty}\xi(\mathbf{z},n)$ 
be the total local time at $\mathbf{z}$ of the infinite path.  

The maximal local time of the walk up to time $n$ is defined as
\begin{equation}
\xi(n):=\max_{\mathbf{z} \in \mathcal{Z}_d}\xi(\mathbf{z},n),
\quad n=1,2,\ldots
\label{loc2}
\end{equation}

Define also the following quantities:
\begin{equation}
\eta(n):=\max_{0\leq k\leq n}\xi(\mathbf{S}_k,\infty), \quad
n=1,2,\ldots
\label{loc3}
\end{equation}

Denote by $\gamma(n)=\gamma(n;d)$ the probability that in the
first $n-1$ steps the $d$-dimensional path does not return to the
origin. Then
\begin{equation}
1=\gamma(1)\geq \gamma(2)\geq ...\geq \gamma(n)\geq...>0.
\label{gain}
\end{equation}
It was proved in \cite{DE50} that

\medskip\noindent
\textbf{Theorem A} (Dvoretzky and Erd\H os \cite{DE50})
\textit{For $d\geq 3$
\begin{equation}
\lim_{n\to\infty}\gamma(n) =\gamma=\gamma(\infty;d)>0,
\label{gam1}
\end{equation}
and
\begin{equation}
{\gamma}<\gamma(n)<{\gamma}+O(n^{1-d/2}).  \label{gam2}
\end{equation}
Consequently
\begin{equation}
\mathbf{P}(\xi(\mathbf{0},n)=0,\,
\xi(\mathbf{0},\infty)>0)=O\left( n^{1-d/2}\right)
\end{equation}
as} $n\to\infty$.

So $\gamma$ is the probability that the $d$-dimensional simple
symmetric random walk never returns to its starting point.

For $d\geq 3$ (see Erd\H os and Taylor \cite{ET60})
$\,\xi(\mathbf{0},\infty)$ has geometric distribution:

\begin{equation}
\mathbf{P}(\xi(\mathbf{0},\infty)=k)={\gamma}(1-{\gamma})^k,\qquad
k=0,1,2,\ldots  \label{geo1}
\end{equation}
Erd\H{o}s and Taylor \cite{ET60} proved the following strong law
for the maximal local time:

\bigskip \noindent \textbf{Theorem B} (Erd\H os and Taylor \cite{ET60})
\textit{For $d\ge 3$
\begin{equation}
\lim_{n\to\infty}\frac{\xi(n)}{\log n}=\lambda \hspace{1cm}
\mathrm{a.s.}, \label{la}
\end{equation}
where}
\begin{equation}
\lambda=\lambda_d=-\frac{1}{\log(1-{\gamma})}.
\end{equation}
We remark that (\ref{la}) is also true if $\xi(n)$ is replaced by
$\eta(n).$

\medskip
Now we present some of our results for the local and
occupation times for simple symmetric random walk in $\zd$, $d\geq
3$. We note however that Theorems E and H  and  the Proposition
below are true for more general symmetric aperiodic random walk in
$\zd$, $d\geq 3$.

\medskip

\bigskip Erd\H{o}s and Taylor \cite{ET60} also investigated the
properties of
\[
Q(k,n):=\#\{\mathbf{z}:\ \mathbf{z}\in\mathcal{Z}_d,\
\xi(\mathbf{z},n)=k\},
\]
i.e. the cardinality of the set of points visited exactly $k$
times in the time interval $[1,n]$.

\medskip\noindent \textbf{Theorem C} (Erd\H os and Taylor \cite{ET60})
\textit{ For $d\geq 3$ and for any} $k=1,2,\ldots$
\begin{equation}
\lim_{n\rightarrow\infty}{\frac{Q(k,n)}{n}}=
\gamma^2(1-\gamma)^{k-1} \hspace{1cm} \mathrm{a.s.}  \label{qkn}
\end{equation}

This was extended in \cite{CsFR06} to a uniform law of large
numbers:

\medskip\noindent \textbf{Theorem D} (Cs\'aki, F\"oldes and R\'ev\'esz
\cite{CsFR06}) \textit{Let $d\geq 3$, and define
\begin{eqnarray}
\mu=\mu(t)&:=&\gamma(1-\gamma)^{t-1},\label{mut}\\
t_n&:=&[\lambda\log n-\lambda B\log\log n],\quad B>2.\label{tn}
\end{eqnarray}
Then we have}
\begin{equation}
\lim_{n\rightarrow\infty}\max_{t\leq
t_n}\left|{\frac{Q(t,n)}{n\gamma\mu(t)}}- 1\right|=0\hspace{1cm}
\mathrm{a.s.}
\end{equation}

\bigskip We introduce the following notations. For
$\mathbf{z}\in \mathcal{Z}_d$ let $T_\mathbf{z}$ be the first hitting
time of $\mathbf{z}$, i.e. $T_%
\mathbf{z}:=\min\{i\geq 1:\mathbf{S}_i=\mathbf{z}\}$ with the
convention that $T_\mathbf{z}=\infty$ if there is no $i$ with
$\mathbf{S}_i=\mathbf{z}$. Let $T=T_\bo$. In general, for a subset
$A$ of $\mathcal{Z}_d$, let $T_A$ denote the first time the random
walk visits $A$, i.e. $T_A:=\min\{i\geq 1:\, \mathbf{S}_i\in A\}
=\min_{\mathbf{z}\in A}T_{\mathbf{z}}$. Let $\mathbf{P}_{
\mathbf{z}}(\cdot)$ denote the probability of the event in the
brackets under the condition that the random walk starts from
$\mathbf{z}\in\mathcal{Z}_d$. We denote $\mathbf{P}(\cdot)=
\mathbf{P}_{\mathbf{0}}(\cdot)$. Define \begg
\gamma_{\mathbf{z}}:=\pe(T_{\mathbf{z}}=\infty).
\endd
Let $B(r)$ be the sphere of radius $r$ centered at the origin, i.e.
\[
B(r):=\{\mathbf{z}\in\mathcal{Z}_d: \Vert \mathbf{z}\Vert=r\},
\]
and $B:=B(1)$ where $\Vert\cdot\Vert$ is the Euclidean norm.

Introduce further
\begin{equation}
\mathbf{p}:=\mathbf{P}_{\mathbf{e}_1}(T_{B}<T).  \label{ap}
\end{equation}

In words, $\mathbf{p}$ is the probability that the random walk, starting
from $\mathbf{e}_1$ (or any other points of $B$),
returns to $B$ before reaching $\mathbf{0}$
(including the case $T_{B}<T=\infty$). It is not
hard to show that
\begg
\mathbf{p}=1-\frac{1}{2d(1-\gamma)}. \label{pe}
\endd

\bigskip
 For a set
$A\subset\mathcal{Z}_d$ the occupation time of $A$ is defined by
\begin{equation}
\Xi(A,n):=\sum_{\mathbf{z}\in A}\xi(\mathbf{z},n).  \label{occup}
\end{equation}

Consider the translates of $A$, i.e. $A+\mathbf{u}=\{\mathbf{z}+\mathbf{u}%
:\, \mathbf{z}\in A\}$ with $\mathbf{u}\in\mathcal{Z}_d$ and
define the maximum occupation time by
\begin{equation}
\Xi^*(A,n):=\max_{\mathbf{u}\in \mathcal{Z}_d}\Xi(A+\mathbf{u},n).
\label{maxop}
\end{equation}
It was shown in {\rm \cite{CsFRRS}}

\bigskip\noindent
\medskip\noindent \textbf{Theorem E} (Cs\'aki, F\"oldes, R\'ev\'esz, Rosen
and Shi {\rm \cite{CsFRRS}})
 \textit {For $d\geq 3$ and for any fixed finite set
$A\subset\mathcal{Z}_d$
\begin{equation}
\lim_{n\to \infty}\frac{\Xi^*(A,n)}{\log n}=
\frac{-1}{\log(1-1/\Lambda_A)}
\qquad\hspace{1cm} \mathrm{%
a.s.},  \label{lim}
\end{equation}
where $\Lambda_A$ is the largest eigenvalue of the $|A|\times |A|$
matrix with elements
$$
G(\mathbf{z}-\mathbf{u}), \qquad
\mathbf{z},\mathbf{u}\in A,
$$
and
$$
G(\mathbf{z})=\sum_{i=0}^\infty
\bp(\mathbf{S}_i=\mathbf{z}),\qquad \mathbf{z}\in\zd
$$
is the Green function of the walk.}

As a major tool for the above result it was proved that

\medskip\noindent \textbf{Proposition A} (Cs\'aki, F\"oldes, R\'ev\'esz,
Rosen and Shi \cite{CsFRRS})
 \begg
\bp(\Xi(A,\infty)>k)=\sum_jh_j\left(\frac{\lambda_j-1}{\lambda_j}\right)^k,
\quad k=0,1,\ldots, \label{XiA}
\endd
\textit{where $\lambda_j$ are the eigenvalues of the matrix $G_A$
and $h_j$ are certain coefficients calculated in terms of the
eigenvectors.}

\medskip
In particular, it was shown in \cite{CsFRRS}

\begin{equation}
\lim_{n\to \infty}\frac{\max_{\kz\in\zd}\Xi(\kz,n)}{\log n}=
\frac{-1}{\log\left(\mathbf{p}+\frac1{2d}\right)}=:\kappa\hspace{1cm}
\mathrm{a.s.}, \label{lim1}
\end{equation}
where
\begin{equation}
\Xi(\kz,n):=\Xi(B+\kz,n),
\label{xiz}
\end{equation}
i.e. the occupation time of the unit sphere centered at $\kz$.
Note that in this notation $\kz$ stands for the center of the unit
sphere not for the one element set $\{\kz\}.$

Furthermore if $A=\{\mathbf{0},\mathbf{z}\}$ is a two-point set,
then the constant  in (\ref{lim}) of Theorem E is
$$
c_A=\frac{-1}{\log(1-1/\Lambda_A)}=
\frac{-1}{\log\left(1-\frac{\gamma}{2-\gamma_{\mathbf{z}}}\right)},
$$
where $\gamma_{\mathbf{z}}$ is the probability that the random
walk, starting from zero, never visits $\mathbf{z}$. It can be seen that 
in this case $c_A<2\lambda$, showing that for large $n$ any point with 
fixed distance from a maximally visited point can not be maximally
visited. This suggests to investigate the behavior of local and
occupation times around frequently visited points.

 Naturally it would be
interesting to investigate  the joint behavior of the local time
and/or  occupation time of two sets in general. However this is a
very difficult proposition. From the  two special cases we
discussed in  {\rm \cite{CsFR1}}  we  mention here the following
one. Consider  the joint behavior of local time of a point and the
occupation time of the unit sphere centered at the point.

Define the set $\mathcal{B}$ on the plane as
\begin{equation}
\mathcal{B}:=\{(x,y):\, y\geq x\geq 0;\,\, -y\log y+x\log
(2dx)+(y-x)\log((y-x)/\mathbf{p})\leq 1\}, \label{defd}
\end{equation}
where $\mathbf{p}$ was defined in (\ref{ap}) and its value in terms of
$\gamma$ is given by (\ref{pe}).

 \medskip\noindent \textbf{Theorem F} (Cs\'aki, F\"oldes and R\'ev\'esz
\cite{CsFR1}) \textit{ Let $d\geq 4$. For each $\varepsilon>0$
with probability  1 there exists an $n_0=n_0(\varepsilon)$ such
that if $n\geq n_0$ then }
\begin{itemize}
\item{\rm (i)} $(\xi(\mathbf{z},n),\Xi(\mathbf{z},n))\in
((1+\varepsilon)\log n) \mathcal{B},\quad
\quad\forall\mathbf{z}\in\mathcal{Z}_d$ \item{\rm (ii)} for any
$(k,\ell)\in ((1-\varepsilon)\log n)\mathcal{B} \cap
\mathcal{Z}_2$ there exists a random $\mathbf{z}\in\mathcal{Z}_2$
for which \end{itemize}
$$
(\xi(\mathbf{z},n),\Xi(\mathbf{z},n))=(k,\ell+1).
$$

\medskip\noindent \textbf{Theorem G} (Cs\'aki, F\"oldes and R\'ev\'esz
\cite{CsFR1})
\textit{  Let $d\geq 3$. For each $\varepsilon>0$ with
probability  1 there exists an $n_0=n_0(\varepsilon)$ such that if
$n\geq n_0$ then}
\begin{itemize}
\item{\rm (i)}
$(\xi(\mathbf{S}_j+\bee_i,\infty),\Xi(\mathbf{S}_j+\bee_i,\infty))\in
((1+\varepsilon)\log n) \mathcal{B},\quad \forall
j=1,2,\ldots,n,\quad \forall i=1,2,\ldots,2d$
\item{\rm (ii)}
for any $(k,\ell)\in ((1-\varepsilon)\log n)\mathcal{B} \cap
\mathcal{Z}_2$ and for arbitrary $i\in\{1,2,\ldots,2d\}$ there
exists a random integer $j=j(k,\ell)\leq n$ for which
\[
(\xi(\mathbf{S}_j+\bee_i,\infty),\Xi(\mathbf{S}_j+\bee_i,\infty))=(k,\ell+1).
\]
\end{itemize}

It follows from these results that if the local time of a point is
close to $\lambda\log n$, then the local times of each of its
neighbors should be asymptotically equal to $\lambda(1-\gamma)\log
n$ which is strictly less than $\lambda\log n$. In \cite{CsFR2} we
investigated whether similar results are true in a wider
neighborhood, i.e. whether the local times of points on a certain
sphere centered at a heavy point are asymptotically determined. As a
positive answer we proved

\medskip\noindent \textbf{Theorem H} (Cs\'aki, F\"oldes and R\'ev\'esz
\cite{CsFR2})
 \textit{  Let $d\geq 5$ and
$k_n=(1-\delta_n)\lambda\log n$. Let $r_n>0$ and $\delta_n>0$ be
selected such that  $\delta_n$ is non-increasing, $r_n$ is
non-decreasing, and for any $c>0$ $r_{[cn]}/r_n<C$ with some $C>0$
and for \begg \beta_n:=r_n^{2d-4}\frac{\log\log n}{\log n}
\label{beta1}
\endd
\begg
\lim_{n\to\infty}\beta_n=0,\qquad \lim_{n\to\infty}\delta_n
r_n^{2d-4}=0.
\endd
Define the random set of points
\begg
{\cal A}_n=\{\mathbf{z}\in Z^d: \xi(\mathbf{z},n)\geq k_n\}.
\endd
Then we have \begg \lim_{n\to\infty}\max_{\mathbf{z}\in {\cal
A}_n} \max_{\mathbf{u}\in
S(r_n)}\left|\frac{\xi(\mathbf{z}+\mathbf{u},n)}
{m_\mathbf{u}\lambda\log n}-1\right|=0\quad{\rm a.s.},
\endd
where}
$$
S(r):=\{\mathbf{u}\in\zd:\, \Vert \mathbf{u}\Vert\leq r\}\quad
{\it and}\quad m_\mathbf{u}:=\be(\xi(\mathbf{u},\infty)\mid
T<\infty)= \frac{(1-\gamma_{\mathbf{u}})^2}{1-\gamma}.
$$

\medskip\noindent \textbf{Theorem I} (Cs\'aki, F\"oldes and R\'ev\'esz
\cite{CsFR2})
\textit{Let $d\geq 3$ and $k_n=(1-\delta_n)\lambda\log n$. Let $r_n>0$
and $\delta_n>0$ be selected such that $\delta_n$ is
non-increasing, $r_n$ is non-decreasing, and for any $c>0$
$r_{[cn]}/r_n<C$ for some $C>0$ and for 
\begg
\beta_n:=r_n^{2d-4}\frac{\log\log n}{\log n} \label{beta2}
\endd
\begg
\lim_{n\to\infty}\beta_n=0,\qquad \lim_{n\to\infty}\delta_n
r_n^{2d-4}=0.
\endd
Define the random set of indices
\begg
{\cal B}_n=\{j\leq n: \xi(\mathbf{S}_j,\infty)\geq k_n\}.
\endd
Then we have} 
\begg \lim_{n\to\infty}\max_{j\in {\cal
B}_n}\max_{\mathbf{u}\in S(r_n)}
\left|\frac{\xi(\mathbf{S}_j+\mathbf{u},\infty)}{m_\mathbf{u}\lambda\log
n}-1\right|=0\quad{\rm a.s.}
\endd

\section{Simple asymmetric random walk on the line}

\renewcommand{\thesection}{\arabic{section}} \setcounter{equation}{0}
\setcounter{theorem}{0} \setcounter{lemma}{0}

\noindent Consider a simple asymmetric random walk on the line
$\{S_n\}_{n=0}^{\infty}$ starting at the origin, i.e.
$S_0:=0$, $S_n:=\sum_{k=1}^{n} X_k$, $n=1,2,\dots$, where $X_k,\,
k=1,2,\dots$ are i.i.d. random variables with distribution
\begg
\bp (X_1=1)=p \,\,\,{\rm and}\,\,\, \bp(X_1= -1)=q\, (=1-p).
\endd
Without restricting  generality we will suppose throughout the paper
that
$$p>q\qquad{\rm and\,\, introduce}\quad h:=\frac{q}{p} .$$

As it is well-known, this random walk is transient, i.e. with
probability one we have $\lim_{n\to\infty}S_n=\infty$. There is a
huge literature on such transient random walk. Some basic results
are given in Feller \cite{F1}, \cite{F2}, Jordan \cite{J}, Spitzer
\cite{S76}, etc., some of these will be given in the next section.

Let ${\cal Z}={\cal Z}_1$, i.e. the set of integers on the line. We
define the local time by
\begin{equation}
\xi(z,n):=\#\{k: \,\,0< k \leq n,\,\,\, S_k=z\},
\quad z\in{\cal Z},\quad n=1,2,\ldots,
\label{lo}
\end{equation}
$$
\xi(z,\infty):=\lim_{n\to\infty}\xi(z,n),\quad
\xi(n):=\max_{z\in{\cal Z}}\xi(z,n),\quad
\eta(n):=\max_{0\leq j\leq n}\xi(S_j,\infty),
$$
and the occupation time of a set $A\subset{\cal Z}$ by
$$
\Xi(A,n):=\#\{k: \,\, 0<k\leq n,\,\,\, S_k\in A\}
=\sum_{z\in A}\xi(z,n),\quad n=1,2,\ldots
$$

Concerning limit theorems for the local time, the analogue of Theorem
A is simple, it will be given as Fact 2 in Section 4. It seems
that no analogue of Theorem B can be found in the literature,
though the following result can be proved from (\ref{zid})-
(\ref{zid2}) just as Theorem B above of Erd\H os and Taylor. It
will be also a trivial consequence of Theorem 3.3 and Theorem 3.4.

\begin{theorem} For the simple asymmetric random walk
\begg \lim_{n\to\infty}\frac{\xi(n)}{\log n}=
\lim_{n\to\infty}\frac{\eta(n)}{\log
n}=\frac{-1}{\log(2q)}=:\lambda_0\quad {\rm a.s.}
\endd
\end{theorem}

We do not know analogues of (\ref{lim}) and (\ref{XiA}), but in
certain particular cases the distribution can be expressed in a
simple form. In the next section we give a version of the joint
distribution of $\xi(z,\infty)$ and $\xi(0,\infty)$   and the
distribution of $\Xi(\{0,z\},\infty)$ (Proposition 4.1).
Furthermore, we present the joint distribution of
$\xi(0,\infty)$ and $\Xi(0,\infty)$ (Proposition 4.2).
Let
$$
\Xi^*(A,n):=\max_{a\in\mathcal{Z}}\Xi(A+a,n).
$$

From (\ref{2point}) in Section 4, similarly to the proof of (\ref{lim})
via (\ref{XiA}) in \cite{CsFRRS}, we will show the following result.
\begin{theorem} For $z>0$ integer
\begg
\lim_{n\to\infty}\frac{\Xi^*(\{0,z\},n)}{\log n}=
\frac{-1}{\log\left(\frac{2q+h^{z/2}}{1+h^{z/2}}\right)} \quad
{\rm a.s.}
\label{xi0z}
\endd
\end{theorem}

Concerning Theorem C, define
\[
\widetilde Q(k,n):=\#\{z\in\mathcal{Z}:\
\xi(z,n)=k\}.
\]
We have from Pitt \cite{Pi} that

\medskip\noindent \textbf{Theorem  K} (Pitt \cite{Pi}) \textit{For}
$k=1,2,\ldots$
$$
\lim_{n\to\infty}\frac{\widetilde
Q(k,n)}{n}=(1-2q)^2(2q)^{k-1}\qquad{\rm
a.s.}
$$

The analogue of Theorem D, i.e. uniform law of large numbers for
$\widetilde Q(k,n)$ remains an open problem.

Next we formulate two theorems which correspond to Theorems F and
G  for the transient walk on the line.

Let $B:=\{-1,1\}$, the one dimensional unit sphere around the
origin. Just like in the  higher dimensional situation we will
denote
\[
\Xi(z,n):=\Xi(B+z,n)=\Xi(\{z-1,z+1\},n)=\xi(z-1,n)+\xi(z+1,n),
\]
i.e. the occupation time of the unit sphere centered at
$z\in\mathcal{ Z}$.

Introduce
\begg g(x,y):=x\log x -y\log
y+(y-x)\log(y-x)-x\log(2p)-y\log q \label{gxy}
\endd
and define the set ${\cal D}$ by
\begin{equation}
{\cal D}:=\{(x,y):y\geq x\geq 0;\,\, g(x,y)\leq 1\}. \label{setd}
\end{equation}

\begin{theorem}
For each $\varepsilon>0$ with probability {\rm 1} there exists an
$n_0=n_0(\varepsilon)$ such that if $n\geq n_0$ then   
\item {\rm(i)} $(\xi(z,n),\Xi(z,n))\in ((1+\varepsilon)\log n)
\mathcal{D},\quad \quad\forall z\in\mathcal{Z}$ 
\item {\rm (ii)} for any $(k,\ell)\in ((1-\varepsilon)\log n)\mathcal{D} 
\cap\mathcal{Z}_2$ there exists a random $z\in\mathcal{Z}$ for which
\[
(\xi(z,n),\Xi(z,n))=(k+1,\ell+2).
\]
\end{theorem}

\begin{theorem} For each $\varepsilon>0$ with probability {\rm 1} there
exists an $%
n_0=n_0(\varepsilon)$ such that if $n\geq n_0$ then

\item  {\rm (i)} $(\xi(S_j,\infty),\Xi(S_j,\infty))\in
((1+\varepsilon)\log n) \mathcal{D},\quad \forall
j=1,2,\ldots,n$

\item  {\rm (ii)} for any $(k,\ell)\in ((1-\varepsilon)\log
n)\mathcal{D}%
\cap \mathcal{Z}_2$ there exists a random integer $j=j(k,\ell)\leq n$

for which
\[
(\xi(S_j,\infty),\Xi(S_j,\infty))=(k+1,\ell+2).
\]
\end{theorem}

From Theorem 3.3 the following consequence is easily obtained.
\begin{corollary}
With probability 1 for all possible sequence of integers $\{z_n\}$ the
set of all possible limit points of
$$
\left(\frac{\xi(z_n,n)}{\log n},\frac{\Xi(z_n,n)}{\log n}\right),
\qquad n\to\infty
$$
is equal to ${\cal D}$.
\end{corollary}

Finally, we state the following analogues of Theorems H and I:

\begin{theorem}
Define the random set of indices \begg {\cal A}_n:=\{u\in {\cal
Z}:\quad \xi(u,n)\geq (1-\delta_n)\lambda_0\log n\}.
\endd

Let $\alpha=\log(1/h),$ and select $c>0$ such that $\alpha c <1.$
If
 $$\lim_{n\to\infty}\delta_n(\log n)^{\alpha c}=0, $$
then we have
\begg \lim_{n\to\infty}\max_{u\in {\cal
A}_n}\max_{|z|\leq c\log\log n} 
\left|\frac{\xi(u+z,n)}{m_z\lambda_0\log n}-1\right|=0\quad{\rm a.s.},
\endd
where \begg m_z :=\left\{
\begin{array}{ll}
& \frac{h^{|z|}}{2q} \ \quad {\rm if} \quad z\neq 0,\\
& 1  \quad\quad \, \, {\rm if} \quad z=0.
\end{array}
\right. \label{mz}
\endd
\end{theorem}

\begin{theorem}
Define the random set of indices
\begg {\cal B}_n:=\{j\leq n:
\xi(S_j,\infty)\geq (1-\delta_n)\lambda_0\log n\}.
\endd
Let $\alpha=\log(1/h),$ and select $c>0$ such that $\alpha c <1.$
If
 $$\lim_{n\to\infty}\delta_n(\log n)^{\alpha c}=0, $$
then we have
\begg \lim_{n\to\infty}\max_{j\in {\cal
B}_n}\max_{|z|\leq c\log\log n} \left|\frac{\xi(S_j+z,\infty)}{
m_z\lambda_0\log n}-1\right|=0\quad{\rm a.s.},
\endd
where $m_z$ is defined in {\rm (\ref{mz})}.
\end{theorem}

\begin{corollary}
Let $A\subset {\cal Z}$ be a fixed set.

{\rm (i)} If  $u_n\in {\cal A}_n$, $n=1,2,\ldots$, then
$$
\lim_{n\to\infty}\frac{\Xi(A+u_n,n)}{\log n}=
\lim_{n\to\infty}\frac{\sum_{x\in A}\xi(x+u_n,n)}{\log n}
=\lambda_0\sum_{x\in A}m_x \quad {\rm a.s.}
$$

{\rm (ii)} If  $j_n\in {\cal B}_n$, $n=1,2,\ldots$, then
$$
\lim_{n\to\infty}\frac{\Xi(A+S_{j_n},\infty)}{\log n}=
\lim_{n\to\infty}\frac{\sum_{x\in A}\xi(x+S_{j_n},\infty)}{\log n}
=\lambda_0\sum_{x\in A}m_x \quad {\rm a.s.}
$$
\end{corollary}

\section{Preliminary facts and results}

\renewcommand{\thesection}{\arabic{section}} \setcounter{equation}{0}
\setcounter{theorem}{0} \setcounter{lemma}{0}

\medskip\noindent
\textbf{Fact 1.} For the probability of no
return  we have (cf. Feller \cite{F1})
$$
\bp(S_i\neq 0,\, i=1,2,\ldots)=1-2q=p-q=:\gamma_0
$$
Let
$$
T_z:=\min\{i\geq 1:\, S_i=z\},\quad T_0=:T,\quad z\in{\cal Z},
$$
the first hitting time of $z$. Denote by $\gamma_0(n)$ the probability
that in the first $n-1$ steps $S_n$ does not return to the origin. Then
just like in (\ref{gain}), we have
\begin{equation}
1=\gamma_0(1)\geq \gamma_0(2)\geq ...\geq \gamma_0(n)\geq...>
\gamma_0>0.
\label{gainq}
\end{equation}

\medskip\noindent
\textbf{Fact 2.} For $T$, the first return time to $0$ we have (cf.
Feller \cite{F1})
\begg \bp(T=2n)={2n\choose n}\frac{1}{2n-1}(pq)^n
 \sim \frac{(4pq)^n}{2\sqrt{\pi}n^{3/2}},\quad
n\to\infty,
\endd
from which one can easily obtain that \begg
\gamma_0(n)-\gamma_0=\bp(n\leq
T<\infty)=O\left(\frac{(4pq)^{n/2}}{n^{3/2}}\right), \quad
n\to\infty. \label{kozel}
\endd
{\bf Remark:} This is the analogue  of Theorem A.

It can be seen furthermore that
\begin{equation}
\bp(T=2n, S_1=1)=\bp(T=2n, S_1=-1)=
\frac12{2n\choose n}\frac{1}{2n-1}(pq)^n=\frac12\bp(T=2n),
\end{equation}
from which one easily obtains
\begin{equation}
\bp(T<\infty, S_1=1)=\bp(T<\infty, S_1=-1)=\frac12 \bp(T<\infty)=q,
\end{equation}
and
\begin{equation}
\bp(n\leq T<\infty, S_1=1)=\bp(n\leq T<\infty, S_1=-1)
=\frac12 \bp(n\leq T<\infty).
\end{equation}

Now introduce
$$q(n):= \bp(T<n, S_1=1)=
\mathbf{P}( T<\infty, S_1=1)-\mathbf{P}(n \leq T<\infty,
 S_1=1)$$
\begg
=\bp(T<n, S_1=-1)=
\mathbf{P} (T<\infty, S_1=-1)-\mathbf{P}(n \leq T<\infty,
S_1=-1).
\label{remm}
\endd
Then we have
\begg
0<q-q(n)=O\left(\frac{(4pq)^{n/2}}{n^{3/2}}\right),\quad
n\to\infty
\label{alp1}
\endd
as well.

Recall the notation $h=q/p (<1).$

\medskip\noindent \textbf{Fact 3.} (see e.g. Feller \cite{F1})
For $z\in {\cal Z}$ we have
\begg \bp(T_z<\infty) =\left\{
\begin{array}{lll}
& h^{-z} \ \quad {\rm if} \quad z<0,\\
& 2q  \quad\quad {\rm if} \quad z=0,\\
& 1   \quad\quad \, \ {\rm if} \quad z>0.
\end{array}
\right. \label{tz}
\endd
\medskip\noindent
\textbf{Fact 4.} (see e.g. Spitzer \cite{S76}, page 10)
For the Green function $G(z)$ we have
 for $z\in {\cal Z}$ : \begg G(z)=\sum_{i=0}^\infty \bp(S_i=z)
=\left\{
\begin{array}{ll}
& \frac{1}{p-q}h^{-z} \quad\, \, {\rm if}\quad z\leq 0, \\
& \frac{1}{p-q}\quad\quad\quad {\rm if} \quad z>0.
\end{array}
\right. \label{gz}
\endd

\begin{lemma}
For $n\geq 1$, $|j|\leq n$ we have
 $$\mathbf{P}(S_n=j)\leq C_1\exp(-C_2 n+C_3j),$$
where the constants $C_i>0$, $i=1,2,3$, depend only on $p.$
\end{lemma}
{\bf Proof.} Clearly  $\mathbf{P}(S_n=j)$ differs from $0$ only if $j$ and
$n$ have the same parity. So we will suppose that in the proof.
\begin{eqnarray*}\mathbf{P}(S_n=j)&=&{ n \choose \frac{n+j}{2}}
p^{\frac{n+j}{2}}q^{\frac{n-j}{2}}
\leq {n \choose [n/2]}(pq)^{n/2}\left(\frac{p}{q}\right)^{j/2}\\
&\leq& C_1(4pq)^{n/2}\left(\frac{p}{q}\right)^{j/2}= C_1\exp(-C_2
n+C_3j),
\end{eqnarray*}
where
$$
C_2=-\frac12\log(4pq),\qquad C_3=\frac12\log\frac{p}{q}.
$$
$\Box$

For the distribution of the local time we have

\medskip\noindent
\textbf{Fact 5.} (cf. Dwass \cite{D})
\begg \bp(\xi(0,\infty)=k)=(2q)^k(1-2q),
\qquad k=0,1,2,\ldots \label{zid}
\endd
For $z>0$ integer
\begg \bp(\xi(z,\infty)=k)=(2q)^{k-1}(1-2q),  \qquad
k=1,2,\ldots \label{zid1}
\endd
 and  \begg \bp(\xi(-z,\infty)=k)=\left\{
\begin{array}{ll}
& 1-h^z \quad \quad \qquad\qquad {\rm if}\quad k=0,\\
& h^z (2q)^{k-1}(1-2q)\quad {\rm if}\quad k=1,2,\ldots
\label{zid2}
\end{array}
\right.
\endd

For the joint distribution of $\xi(z,\infty)$ and $\xi(0,\infty)$
we have
\begin{propo} For $z>0$, $k\geq 0$ integers we have

\begg \be\left(e^{v\xi(z,\infty)},\,\xi(0,\infty)=k\right)
=(1-2q)(2q)^k\varphi^k(v)\psi(v),
\label{joint}
\endd
\begg \be\left(e^{v\xi(-z,\infty)},\,\xi(0,\infty)=k\right)
=(1-2q)(2q)^k\varphi^k(v), \quad \label{ojoint}
\endd
for
$$
v<-\log\left(1-\frac{1-2q}{1-h^z}\right),
$$
where
\begg
\varphi(v):=\frac{1-\frac{4q^2-h^z}{2q(1-2q)}\left(e^v-1\right)}
{1-\frac{2q-h^z}{(1-2q)}\left(e^v-1\right)},
\endd
\begg \psi(v):=\frac{e^v}{1-\frac{2q-h^z}{(1-2q)}\left(e^v-1\right)}.
\endd

\noindent
Moreover,
\begg \bp(\Xi(\{0,z\},\infty)=k)=\frac{1-2q}{2h^{z/2}}
\left(\left(\frac{2q+h^{z/2}}{1+h^{z/2}}\right)^k-
\left(\frac{2q-h^{z/2}}{1-h^{z/2}}\right)^k\right), \label{2point}
\endd
$k=1,2,\ldots$

\begg \bp(\Xi(\{0,-z\},\infty)=k)=\frac{1-2q}{2}
\left(\left(\frac{2q+h^{z/2}}{1+h^{z/2}}\right)^k+
\left(\frac{2q-h^{z/2}}{1-h^{z/2}}\right)^k\right),
\label{2point2}
\endd
$k=0,1,2,\ldots$
\end{propo}

\noindent{\bf Proof.}
Let us recall the gambler ruin (cf. Feller \cite{F1} or Jordan
\cite{J}): for $0\leq a<b<c$
\begg \mathbf{P}_b(T_a<T_c)= 1-\frac{1-h^{b-a}}{1-h^{c-a}}.
\label{garu}\endd
Let $z>0$ be an integer. Then by (\ref{garu})

$$s_z:=\bp(T_z<T)=p\,\bp_1(T_z<T)=p\frac{1-h}{1-h^z}=:P_z.
$$
On the other hand,
$$s_{-z}=\bp(T_{-z}<T)=\bp_z(T<T_z)=q\frac{h^{z-1}-h^z}{1-h^z}=h^zP_z.$$
Similarly, a simple calculation shows
$$q_z:=\bp(T<T_z)=\bp_z(T_z<T)=1-P_z=: Q_z,  \qquad {\rm and} \qquad
q_{-z}=\bp(T<T_{-z})=q_z.
$$

Let $Z(A)$ be the number of visits in the set
$A\subset{\cal Z}$ in the first excursion away from $0$. In particular, 
for the one point set $\{z\}$
$$ Z(\{z\})=\xi(z,T).$$
Note that $T=\infty$ is possible.

\medskip\noindent
\textbf{Fact 6.} (Baron and Rukhin \cite{BR99})
For $z>0$ integer
\begg \bp(Z(\{z\})=j,\,
T<\infty)=\bp(Z(\{-z\})=j,\, T<\infty)=\left\{
\begin{array}{ll}
& Q_z\quad\quad \quad\quad{\rm if} \quad j=0, \\
& h^z P_z^2 Q_z^{j-1} \quad{\rm if} \quad j=1,2,\ldots \label{d1}
\end{array}
\right.
\endd
\begg \bp(Z(\{z\})=j,\, T=\infty)=(1-2q)P_zQ_z^{j-1},\qquad
j=1,2,\ldots \label{d2}
\endd
\begg \bp(Z(\{-z\})=0,\, T=\infty)=(1-2q).\qquad \label{d3}
\endd

It can be seen furthermore that
$$\be(\xi(z,T),\, T<\infty)=\be(Z(\{z\}),\, T<\infty)=\be(Z(\{-z\}),\,
T<\infty)=h^z,$$
hence
$$\be(\xi(z,T)|T<\infty)=\be(Z(\{z\})|T<\infty)=\be(Z(\{-z\})|T<\infty)
=m_z, $$
where $m_z$ is defined by (\ref{mz}).

Now (\ref{joint}) can be calculated from (\ref{d1}) and (\ref{d2})
by using that
$$
\be\left(e^{v\xi(z,\infty)},\xi(0,\infty)=k\right)=
\left(\be(e^{vZ(\{z\})},\, T<\infty)\right)^k \be(e^{vZ(\{z\})},\,
T=\infty).$$
It is easy to see that
$$\be(e^{vZ(\{z\})},\, T<\infty)=Q_z+\frac{h^z P^2_z
e^v}{1-Q_ze^v}=2q\varphi(v),
$$
and
$$\be(e^{vZ(\{z\})},\,
T=\infty)=(1-2q)\frac{e^vP_z}{1-Q_ze^v}=(1-2q)\psi(v).$$
To ease the computation we remark that
$$h^zP^2_Z-Q^2_z=\frac{h^z-4q^2}{1-h^z.}$$
 Similarly we get (\ref{ojoint})   from  (\ref{d1})
and (\ref{d3}).

From (\ref{joint}) - (\ref{ojoint}), substituting $w=e^v$, one can
find
\begg
\be(w^{\xi(0,\infty)+\xi(z,\infty)})= \frac{(1-2q)P_z w}{1-2Q_z
w+(Q_z^2-h^z P_z^2)w^2},\endd
and
\begg
\be(w^{\xi(0,\infty)+\xi(-z,\infty)})= \frac{(1-2q)(1-Q_z
w)}{1-2Q_z w+(Q_z^2-h^z P_z^2)w^2},\label{joge}
\endd
consequently (\ref{2point}) and (\ref{2point2}) can be obtained by
expanding the right-hand side into powers of $w$. $\Box$

The next result concerns the joint distribution of the local time of
the origin and the occupation time of the  unit sphere $B$.

\begin{propo} For $K=0,1,\ldots,\, L=K+1,K+2,\ldots$ we have
\begg \bp(\Xi(0,\infty)=L, \xi(0,\infty)=K)= {L-1\choose K}(2p)^K
q^{L-1}p(1-2q)=:p(L,K).
\label{d5}
\endd
\begg \bp(\Xi(0,n)=L, \xi(0,n)=K)\leq p(L,K),
\label{d15}
\endd

\begg \bp(\Xi(0,\infty)=L)=(q+2pq)^{L-1}p(1-2q),\quad
L=1,2,\ldots \label{Xid} \endd
Moreover, for the occupation time of the set $A_0:=\{-1,0,1\}$ we have
\begg
\pe(\Xi(A_0,\infty)=\ell)=p(1-2q)\left(\frac{q}{2}\right)^{\ell-1}
\frac{(1+\beta)^\ell-(1-\beta)^\ell}{2\beta},\qquad \ell=1,2,\ldots,
\label{d16}
\endd
where $$\beta=\sqrt{1+\frac{8p}{q}}.$$
\end{propo}
{\bf Proof.} The probability that the walk goes
to $B$ and returns to $0$ immediately, is $2pq.$ However, before
returning, the walk can make one or more outward excursions. If
the walk starts outward from $B$, it returns with probability $q,$
independently whether it starts from $1$ or $-1.$ Altogether we
need $K$ trips from zero to $B$ and back to $0$ and $L-K-1$
outward excursions. Every such arrangement, independently of the
order of occurrences, has probability  $(2pq)^K q^{L-K-1}.$ The
number of ways how we can order the outward excursions and the
inward ones is ${L-1\choose K}.$  After the $K$-th return to zero
the walk must go to $1$ (some of the outward excursion might
happen  at this time) and then to infinity which has probability
$p(1-2q)$, proving our first statement. Now (\ref{d15}) follows
from (\ref{d5}) and (\ref{alp1}).

With a similar argument one can show that if we start the walk at
$-1$, then  for $K=1,2,\ldots, L=K,K+1,\ldots $
\begg
p_{-1}(L,K):=
\mathbf{P}_{-1} (\Xi(0,\infty)=L, \xi(0,\infty)=K)={L\choose K}(2p)^{K-1}
q^{L-1}p^2(1-2q) \label {d6}
\endd
(which does not include the very first visit at $-1$).

Similarly for $K=1,2,\ldots, L=K,K+1,\ldots$
\begg
p_{1}(L,K):=
\mathbf{P}_{1} (\Xi(0,\infty)=L, \xi(0,\infty)=K)={L\choose K}(2p)^{K-1}
q^{L}p(1-2q)
\label{d7}
\endd
(which does not include the very first visit at $1$),
and for $L=0,1,\ldots$
\begg
p_{1}(L,0):= \mathbf{P}_{1} (\Xi(0,\infty)=L, \xi(0,\infty)=0)=
q^L(1-2q).
\endd

Finally, we get (\ref{d16}) from
$$
\be(w^{\Xi(A_0,\infty)})=\be(w^{(\Xi(0,\infty)+\xi(0,\infty))})=
\frac{p(1-2q)w}{1-qw-2pqw^2}
$$
which follows from
$$
\be(e^{vZ(B)}, T<\infty)=\frac{2pq e^v}{1-qe^v}
$$
and
$$
\be(e^{vZ(B)}, T=\infty)=p(1-2q)\frac{e^v}{1-qe^v},
$$
similarly as we got (\ref{2point}) in Proposition  4.1. $\Box$

We will need some basic observations about the reversed walk. By
the reversed path of $(S_0,S_1,\ldots,S_j)$ we mean the path $(0,
S_{j-1}-S_j, S_{j-2}-S_j,\ldots, S_0-S_j)=: (S_0^*,S_1^*,\ldots,
S_j^*)$, e.g. $X_i^*=S_{j-i}-S_{j-i+1}$ for $1\leq i\leq j.$ Then
of course $X_i^*, \quad 1\leq i \leq j$  are i.i.d. random
variables with \begg \bp (X^*_1=1)=q \,\,\,{\rm and}\,\,\, P
(X^*_1= -1)=p=1-q,  \quad  p>q.
\endd
If needed this can be extended to an infinite path $(S_0^*,S_1^*,\ldots,
S_n^*,\ldots)$, i.e. $\{S^*_n\}_{n=0}^{\infty}$ is a walk starting
at the origin $S^*_0=0$, with  $S^*_n=\sum_{k=1}^{n} X^*_k$,
$n=1,2,\ldots$, where $X^*_k,\, $ are defined above for $k\leq j$
and for $k>j$ they are an arbitrary sequence of i.i.d. random variables
(also independent from the previous ones) with the above distribution.
Then we clearly have for all $z$

\medskip\noindent \textbf{Fact 7.}
\begg \mathbf{P}(\xi^*(-z,\infty)=k)=\mathbf{P}(\xi(z,\infty)=k).
\endd

Consequently, under the conditions of Proposition 4.1 we have for
$z=1,2,\ldots$,   $k=0,1,2,\ldots$  and $v<-\log Q_z$

\begg \be\left(e^{v\xi^*(-z,\infty)},\,\xi^*(0,\infty)=k\right)
=(2q)^k(1-2q)\varphi^k(v)\psi(v), \quad  \label{*joint}
\endd
\begg \be\left(e^{v\xi^*(z,\infty)},\,\xi^*(0,\infty)=k\right)
=(2q)^k(1-2q)\varphi^k(v). \quad \label{*ojoint}
\endd

\section{Proof of Theorem 3.2}

\renewcommand{\thesection}{\arabic{section}} \setcounter{equation}{0}
\setcounter{theorem}{0} \setcounter{lemma}{0}

\begin{lemma}
Let
$$
\theta=-\log\frac{2q+h^{z/2}}{1+h^{z/2}}.
$$
There exist $u_0>0$, $c_1>0$, $c_2>0$ such that for all $u\geq u_0$,
$n\geq u^2$ we have
\begg
c_1e^{-\theta u}\leq
\bp(\Xi(\{0,z\},n)\geq u)\leq\bp(\Xi(\{0,z\},\infty)\geq u)\leq
c_2e^{-\theta u}.
\label{inequ}
\endd
\end{lemma}

\noindent \textbf{Proof.} The second inequality in (\ref{inequ}) is
obvious, and since from (\ref{2point}) we can see
$$
\bp(\Xi(\{0,z\},\infty)\geq u)\sim ce^{-\theta u},\qquad u\to\infty,
$$
with some constant $c$, we have also the third inequality in
(\ref{inequ}). To show the first inequality in (\ref{inequ}), we note
that
$$
\bp(\Xi(\{0,z\},\infty)\geq u)\leq \bp(\Xi(\{0,z\},n)\geq u)
+\bp(\cup_{k=n+1}^{\infty}\{S_k=0\})+\bp(\cup_{k=n+1}^{\infty}\{S_k=z\}).
$$
By Lemma 4.1 for $j=0,z$
$$
\bp(\cup_{k=n+1}^{\infty}\{S_k=j\})\leq C_1 e^{-Cn},
$$
with some $C >0, \,C_1>0 $. Thus, by (\ref{2point}) we have
$$
\bp(\Xi(\{0,z\},n)\geq u)\geq\bp(\Xi(\{0,z\},\infty)\geq u)
-2C_1 e^{-Cn}\geq ce^{-\theta u}-2C_1 e^{-Cu^2}\geq C_2e^{-\theta u}
$$
for large enough $u$ with some $C_2>0$. $\Box$

Now we turn to the proof of Theorem 3.2. First we prove an upper
bound in (\ref{xi0z}). Since
$$
\sum_{a=-\infty}^z\xi(a,\infty)<\infty\qquad{\rm and}
\qquad \sum_{a=n+1}^\infty\xi(a,n)=0
$$
almost surely, it suffices to show the upper bound with
$\Xi^*(\{0,z\},n)$ replaced by
$$
\widetilde\Xi^*(\{0,z\},n):=\max_{0\leq a\leq n}\Xi(\{a,a+z\},n).
$$
Since $\bp(\Xi(\{a,a+z\},n)>u)\leq \bp(\Xi(\{0,z\},n)>u-1)$  for
$a=1,2,\ldots,$ we have by Lemma 5.1 for large enough $n$
$$
\bp\left(\widetilde\Xi^*(\{0,z\},n)\geq
\frac{1+\varepsilon}{\theta}\log n\right)\leq (n+1)
\bp\left(\Xi(\{0,z\},n)\geq \frac{1+\varepsilon}{\theta}\log
n -1\right) $$
$$ \leq (n+1) \bp\left(\Xi(\{0,z\},\infty)\geq
\frac{1+\varepsilon}{\theta}\log n-1\right)\leq \frac{c}{n^{\varepsilon}}.
$$
Applying this for the subsequence $n_k=k^{2/\varepsilon}$, using
Borel-Cantelli lemma and the monotonicity of $\Xi^*(\{0,z\})$,
$\varepsilon$ being arbitrary, we obtain
$$
\limsup_{n\to\infty}\frac{\widetilde\Xi^*(\{0,z\},n)}{\log n}
\leq \frac{-1}{\log\left(\frac{2q+h^{z/2}}{1+h^{z/2}}\right)} \quad
{\rm a.s.}
$$
implying the upper bound in (\ref{xi0z}).

Now we show the lower bound in (\ref{xi0z}). Let $k(n)=[\log n]^3$,
$N_n=[n/k(n)]$, $t_{i,n}=ik(n)$, $i=0,1,\ldots, N_n-1$. We have
$$
\Xi^*(\{0,z\},n) \geq \max_{0\leq i\leq N_n-1}Z_i,
$$
where
$$
Z_i:=\Xi(\{S_{t_{i,n}},S_{t_{i,n}}+z\},t_{i+1,n})-
\Xi(\{S_{t_{i,n}},S_{t_{i,n}}+z\},t_{i,n}).
$$
$Z_i,\, i=0,1,\ldots N_n-1$ are i.i.d. random variables, distributed as
$\Xi(\{0,z\},k(n))$, so we get by Lemma 5.1
$$
\bp\left(\Xi^*(\{0,z\},n)\leq \frac{1-\varepsilon}{\theta}\log n\right)
\leq \bp\left(\max_{0\leq i\leq N_n-1}Z_i\leq
\frac{1-\varepsilon}{\theta}\log n\right)
$$
$$
\leq \left(1-\bp\left(\Xi(\{0,z\},k(n))\geq
\frac{1-\varepsilon}{\theta}\log n\right)\right)^{N_n}
\leq\left(1-c_1e^{-(1-\varepsilon)\log n}\right)^{N_n} \leq
e^{-c_1n^{\varepsilon}/k(n)}.
$$
The lower bound in (\ref{xi0z}) follows by Borel-Cantelli lemma.
$\Box$

\section{Proof of Theorems 3.3 and 3.4}

\renewcommand{\thesection}{\arabic{section}} \setcounter{equation}{0}
\setcounter{theorem}{0} \setcounter{lemma}{0}

From (\ref{d5}) (\ref{d6}), (\ref{d7}) and Stirling formula we conclude
the following limit relations.
\begin{lemma} For $y\geq x\geq 0$ we have
$$
\lim_{n\to\infty}\frac{\log p([y\log n]+1,[x\log n])}{\log n}=
\lim_{n\to\infty}\frac{\log p_1([y\log n],[x\log n])}{\log n}
$$

\begg
=\lim_{n\to\infty}\frac{\log p_{-1}([y\log n]+1,[x\log n]+1)}{\log n}
=-g(x,y),
\label{beq1}
\endd
where $g(x,y)$ is defined by {\rm (\ref{gxy})}.
\end{lemma}

Consequently, the probability 
$\pe(\Xi(0,\infty)=[y\log n]+1,\,\xi(0,\infty)=[x \log n])$ is of
order $1/n$, if $(x,y)$ satisfies the basic equation

\begg g(x,y)=1,\qquad y\geq x\geq 0. \label{bas1}
\endd

The following lemma describes the main
properties of the boundary of the set ${\cal D}.$

\newpage

\begin{lemma}
\item {\rm (i)} For the maximum value of $x,y$, satisfying {\rm
(\ref{bas1})}, we have
\begin{eqnarray}
x_{\max}&=&{\frac{-1}{\log (2q)}=\lambda_0,}\label{xmax} \\
y_{\max}&=&{\frac{-1}{\log(q(1+2p))}}=:\kappa_0.
\end{eqnarray}
\item{\rm (ii)} If $x=x_{\max}=\frac{-1}{\log (2q)}$, then
$y=\frac{-1}{p\log (2q)}$. If $y=y_{\max}=\kappa_0$, then $x=(2\kappa_0
p)/(2p+1)$. If $x=0$, then $y=-1/\log q$.

\item {\rm (iii)} For given $x$, the equation {\rm (\ref{bas1})}
has one solution in $y$ for $0\leq x< -1/\log(2pq)$ and for
$x=\lambda_0$, and two solutions in $y$ for $\, -1/\log(2pq)\leq
x<\lambda_0$.
\end{lemma}

\noindent \textbf{Proof.} (i) First consider $x$ as a function
of $y$ satisfying (\ref{bas1}). We seek the maximum, where the
derivative $x^{\prime}(y)=0$. Differentiating (\ref{bas1}) and
putting $x^{\prime}=0$, a simple calculation leads to
\[
-\log y+\log(y-x)-\log q=0,
\]
i.e.
\[
y=x/p.
\]
It can be seen that this is the value of $y$ when $x$ takes its
maximum. Substituting this into (\ref{bas1}), we get
\[
x_{\max}=\frac{-1}{\log (2q)},
\]
verifying (\ref{xmax}).

Next consider $y$ as a function of $x$ and maximize $y$ subject to
(\ref{bas1}). Again, differentiating (\ref{bas1}) with respect to
$x$ and putting $y^{\prime}=0$, we get
\[
- \log(y-x)+\log x-\log(2p)=0
\]
from which $x= (2py)/(1+2p)$. Substituting in (\ref{bas1}) we get
$y_{\max}=\kappa_0$.

This completes the proof of Lemma 6.2(i) and the first two
statements in Lemma 6.2(ii). A simple  calculation shows that if
$x=0$ then $y=-1/\log q$.

Now we turn to the proof of Lemma 6.2(iii). For given $0\leq x\leq
\lambda_0$ consider $g(x,y)$ as a function of $y$. We have
\[
{\frac{\partial g}{\partial y}}=\log{\frac{y-x}{qy}}
\]
and this is equal to zero if $y=x/p$. It is easy to see that $g$
takes a minimum here and is decreasing if $y<x/p$ and increasing
if $y>x/p$. Moreover,
\[
{\frac{\partial^2 g}{\partial
y^2}}={\frac{1}{y-x}}-{\frac{1}{y}}>0,
\]
hence $g$ is convex from below. We have for $0<x<\lambda_0,$ that
this minimum is
\[
g\left(x,\frac{x}{p}\right)=x\log(1/(2q))=\frac{x}{\lambda_0}<1,
\]
and
\[
g(x,x)=-x\log(2pq)\left\{
\begin{array}{ll}
& < 1\,\, \mathrm{if}\,\, x<-1/\log(2pq), \\
& = 1\,\, \mathrm{if}\,\, x=-1/\log(2pq), \\
& > 1\,\, \mathrm{if}\,\, x>-1/\log(2pq).
\end{array}
\right.
\]
This shows that equation (\ref{bas1}) has one solution if $0\leq
x<-1/\log(2pq)$ and two solutions if $-1/\log(2pq)\leq
x<\lambda_0$.

For $x=\lambda_0$, it can be seen that $y=\lambda_0 /p$ is the only
solution of $g(x,y)=1$.

The proof of Lemma 6.2 is complete. $\Box$

\medskip\noindent
{\bf Proof of Theorem 3.4(i).}
Obviously, for $z=1,2,\ldots$ we have for $K=1,2,\ldots,\,
L=K+1,K+2,\ldots$
$$
\bp(\Xi(z,\infty)=L,\, \xi(z,\infty)=K)=
\bp_{-1}(\Xi(0,\infty)=L-1,\, \xi(0,\infty)=K)
$$
and for $K=0,1,\ldots,\, L=K+1,K+2,\ldots$
$$
\bp(\Xi(-z,\infty)=L,\, \xi(-z,\infty)=K)\leq
\bp_1(\Xi(0,\infty)=L-1,\, \xi(0,\infty)=K).
$$

Hence for $(k,\ell)\notin ((1+\ep)\log n) {\cal D}$ and $z\in{\cal Z}$,
as $g(cx,cy)=cg(x,y)$ for any $c>0$,  we have by (\ref{d5})--(\ref{d7})
and Lemma 6.1

\begg
\pe(\xi(z,\infty)=k,\Xi(z,\infty)=\ell)\leq\frac{c}{n^{1+\ep}}.
\label{com1}
\endd
Consequently, by  Fact 5, (\ref{Xid}),  (\ref{com1}) we have

\[
\mathbf{P}(\xi(S_j,\infty),\, \Xi({S _j}, \infty))\not\in
((1+\varepsilon)\log n)\mathcal{D})
\]
\[
\leq \sum_{{(k,\ell)\not\in ((1+\varepsilon)\log n){\cal D} \atop
k\leq (1+\varepsilon)\lambda_0\log n} \atop \ell\leq
(1+\varepsilon)\kappa_0\log n}
\mathbf{P}(\xi(S_j,\infty)=k,\Xi(S_j,\infty)=\ell)
\]
\[
+\sum_{k>(1+\varepsilon)\lambda_0\log n}
\mathbf{P}(\xi(S_j,\infty)=k)
+\sum_{\ell>(1+\varepsilon)\kappa_0\log n}
\mathbf{P}(\Xi(S_j,\infty)=\ell)
\]
\[
\leq\frac{c\log^2 n}{n^{1+\varepsilon}}+
\sum_{k>(1+\varepsilon)\lambda_0\log n}c(2q)^k +
\sum_{\ell>(1+\varepsilon)\kappa_0\log n}c(q+2pq)^\ell
\leq\frac{c}{n^{1+\varepsilon/2}}.
\]
where in the above computation $c$ is an unimportant constant. We
continue denoting such constants by the same letter $c$, the value of
which might change from line to line. Selecting
$n_r=r^{4/\varepsilon},$ we have
\begg \mathbf{P}(\cup_{j\leq
n_{r+1}} \{(\xi(S_j,\infty), \Xi(S_j,\infty))\not\in
((1+\varepsilon)\log n_r)\mathcal{D}\})\leq {\frac{c}{
n_r^{\varepsilon/2}}}=\frac{c}{r^2}.
\endd

This combined with  the Borel-Cantelli lemma shows that with
probability {\rm 1} for all large $r$ and $j\leq n_{r+1}$
$$
(\xi(S_j,\infty),\Xi(S_j,\infty))\in ((1+\varepsilon)\log
n_r)\mathcal{D}.
$$
It follows that with probability 1 there exists an $n_0$ such that if 
$n\geq n_0$ then
$$
(\xi(S_j,\infty),\Xi(S_j,\infty))\in ((1+\varepsilon)\log
n)\mathcal{D}
$$
for all $j\leq n$.

This proves (i) of Theorem 3.4. $\Box$

\bigskip\noindent
{\bf Proof of Theorem 3.3(i).}
The proof is similar to that of Theorem 1.1(i)  in \cite{CsFR1}.

Define the following events for $j\leq n$:
\begin{eqnarray}
B(j,n)&:=& \{(\xi(S_j,n),\Xi(S_j,n))
\notin ((1+\varepsilon)\log n)\mathcal{D}\},\\
B^*(j,n)&:=&
\{(\xi(S_j,j),\Xi(S_j,j))\notin ((1+\varepsilon)\log n)\mathcal{D}\},\\
C(j,n)&:=&\{S_m\neq S_j, m=j+1,\ldots,n
\, \},\\
D(j,n)&:=& \{\Xi(S_j,\infty)>\Xi(S_j,n)\}.
\end{eqnarray}

Considering the reverse random walk starting from $S_j$, i.e.
$S_r'=S_{j-r}-S_{j}$, $r=0,1,\ldots, j$, we remark
$\xi(S_j,j)=\xi'(0,j)$, if $S_j=0$, $\xi(S_j,j)=\xi'(0,j)+1$, if
$S_j\neq 0$, $\Xi(S_j,j)=\Xi'(0,j),$ where $\Xi'$ is the occupation time
of the unit sphere of the random walk $S'$.

From this we can follow the proof of Theorem 1.1(i) in
\cite{CsFR1}, using (\ref{d15}) and (\ref{beq1}) instead of
(2.18) and (3.1) in \cite{CsFR1} and applying Theorem 3.4(i)
instead of Theorem 1.2(i) in \cite{CsFR1}. $\Box$

\medskip\noindent
{\bf Proof of Theorems 3.3(ii) and 3.4(ii)}.
We will say that $S_i$ is new if
\begg \max_{0\leq m<i} S_m<S_i.
\endd

\begin{lemma}
Let $\nu_n$ denote the number of new points up to time n. Then
$$\lim_{n\to \infty} \frac{\nu_n}{n}= 1-2q
\qquad{\rm a.s.}$$
\end{lemma}

\noindent
{\bf Proof}: Let
\[
Z_i:=\left\{
\begin{array}{ll}
1 & \, {\rm if}\,\, S_i \,\, {\rm is}\,\, {\rm new}
\\
0 & \, \mathrm{otherwise.}
\end{array}
\right.
\]
Then $\nu_n=\sum_{i=1}^n Z_i$.

\[ \mathbf{E}(\nu_n^2)=\mathbf{E}\left(\sum_{j=1}^n\sum_{i=1}^n
Z_jZ_i\right)= \mathbf{E}\left(\sum_{j=1}^n Z_j\right)+
2\mathbf{E}\left(\sum_{j=1}^n \sum_{i=1}^{j-1}Z_jZ_i\right)
\]
\[
\leq n
+2\sum_{j=1}^n\sum_{i=1}^{j-1}\mathbf{P}(Z_i=1)\mathbf{P}(Z_{j-i}=1).
\]

Considering the reverse random walk from $S_i$ to $S_0=0$, we see
that the event $\{Z_i=1\}$ is equivalent to the event that this
reversed random walk starting from $0$ does not return to $0$ in
time $i.$  We remark that for the reversed walk the probability of
stepping to the right is $q$ and stepping to left is $p.$
Using (\ref{kozel}) and observing  that it remains true  for the
reversed random walk as well, we get
\[
\mathbf{P}(Z_i=1)=\gamma_0+O((4pq)^{i/2})=1-2q+O((4pq)^{i/2}).
\]

Hence
\[
\mathbf{E}(\nu_n^2)
\leq
n+2\sum_{j=1}^n\sum_{i=1}^{j-1}\left(1-2q+O(4pq)^{i/2}\right)\left((1-2q)
+O(4pq)^{(j-i)/2}\right)
\]
\[
=n(n-1)(1-2q)^2+O(n).
\]
As $\be(\nu_n)=n(1-2q)+O(1)$, we have
 \[ Var(\nu_n)=O(n).
\]
By Chebyshev's inequality we get that
\[ P|\nu_n-n(1-2q)|>\ep n)\leq O\left(\frac{1}{n}\right).  \]
Considering  the subsequence $n_k=k^2$ and using the Borel-Cantelli
lemma and monotonicity of $\nu_n$, we obtain  the lemma.
$\Box$

To show Theorem 3.4(ii), let $\{a_n\}$ and $\{b_n\}$ ($a_n\log
n\ll b_n\ll n$ ) be two sequences to be chosen later. Define
\[
\theta_1=\min\{i>b_n: S_i\, \, {\rm is}\, \, {\rm new}\},
\]
\[
\theta_m=\min\{i>\theta_{m-1}+b_n: S_i\, \, {\rm is} \, \, {\rm new}%
\},\quad m=2,3,\dots
\]
and let $\nu_n^{\prime}$ be the number of $\theta_m$ points up to
time $n-b_n$. Obviously $\nu_n^{\prime}(b_n+1)\geq \nu_n$, hence
$\nu_n^{\prime}\geq\nu_n/(b_n+1)$ and it follows from Lemma 6.3
that for $c<1-2q$, we have with probability 1 that
$\nu_n^{\prime}>u_n:=cn/(b_n+1)$ except for finitely many $n$.

Recall that $B=\{-1,1\}$ denotes the unit sphere around $0$. Let
$$
\rho_0^i=0, \qquad \rho_h^i=\min\{j>\rho_{h-1}^i:
S_{\theta_i+j}\in S_{\theta_i}+B\},\quad h=1,2,\ldots,
$$
i.e. $\rho_h^i,\, h=1,2,\ldots$ are the times when the random walk
visits the unit sphere around $S_{\theta_i}$.

For a fixed pair of integers $(k,\ell)$  define the following
events:
\begin{eqnarray*}
A_i&:=&A_i(k,\ell)=\{\xi(S_{\theta_i}+1,\theta_i+\rho_{\ell+1}^i)=k+1,
\Xi(S_{\theta_i}+1, \theta_i+\rho_{\ell+1}^i)=\ell+2,\\
&&\rho_h^i-\rho_{h-1}^i< a_n,\, h=1,\ldots,\ell+1,\, S_j\not\in
S_{\theta_i}+B,\,
j=\theta_i+\rho_{\ell+1}^i+1,\ldots,\theta_i+b_n\},\\
B_i&:=&B_i(k,\ell)= \{S_j\not\in S_{\theta_i}+B,\,\,
j>\theta_i+b_n\},\\
C_n&:=&C_n(k,\ell)=A_1B_1+\overline{A_1}A_2B_2+\overline{A_1}\,\,
\overline{A_2}A_3B_3+\ldots
+\overline{A_1}\ldots\overline{A_{u_n-1}}A_{u_n}B_{u_n},
\end{eqnarray*}
where $\overline A$ denotes the complement of $A$.

Then we have $\mathbf{P}(A_i)=\mathbf{P}(A_1)$ and
$\mathbf{P}(A_iB_i)=\mathbf{P }(A_1B_1)$, $i=2,3,\dots$ and
\[
\mathbf{P}(C_n)=\mathbf{P}(A_1B_1)\sum_{j=0}^{u_n-1}(1-\mathbf{P}(A_1))^j=
\frac{\mathbf{P}(A_1B_1)}{\mathbf{P}(A_1)}(1-(1-\mathbf{P}(A_1))^{u_n}),
\]
\[
\mathbf{P}(\overline{C_n})\leq
1-\frac{\mathbf{P}(A_1B_1)}{\mathbf{P}(A_1)}%
+e^{-u_n\mathbf{P}(A_1)}.
\]
$A_1B_1$ is the event that starting from the new point
$S_{\theta_1}$, the random walk visits $S_{\theta_1}+1$ exactly
$k+1$ times, while it visits the unit sphere
around this point exactly $\ell+2$ times (including the initial visit at
$S_{\theta_1}$) and all the time intervals between consecutive visits 
are less than $a_n$. Since the return to the sphere via its center takes 
only 2 steps, we have to control only the returns from outside. 
Similarly to (\ref{d6}), one can see
\[
\mathbf{P}(A_1B_1)={{\ell+1}\choose{k+1}}
(q(a_n))^{%
\ell-k}(2pq)^{k}p^2\left(1-2q\right),
\]
\[
={{\ell+1}\choose{k+1}}
(q+O((4pq)^{a_n}))^{%
\ell-k}(2pq)^{k}p^2\left(1-2q\right)
\]
and
\[
\mathbf{P}(A_1)= {{\ell+1} \choose {k+1}}(2pq)^{k}
(q+O((4pq)^{a_n}))^{\ell-k}p^2\left(1-2q+O((4pq)^{b_n-\ell a_n})\right),
\]
where $q(n)$ is given in (\ref{remm}).

Using $a_n\log n\ll b_n\ll n$
\[
\frac{\mathbf{P}(A_1B_1)}{\mathbf{P}(A_1)}=1+O\left((4pq)^{c_1b_n/2}\right),
\]
 for some  $c_1>0$ depending only on $p,$ hence
\[
\mathbf{P}(\overline{C_n})\leq  e^{-c_1 b_n}
+e^{-cn\mathbf{P}(A_1)/b_n}.
\]

For fixed $\varepsilon>0$ introduce the notation
${\cal G}_n=((1-\varepsilon)\log n){\cal D}\cap{\cal Z}_2$.
Choosing $b_n=n^{\delta/2}$, $a_n=n^{\delta/4}$, we
can prove using Stirling formula that for $(k,\ell)\in {\cal G}_n$
 \[
\mathbf{P}(A_1)\geq
\frac1{n^{1-\delta}}
\]
for some $\delta>0$. Since the cardinality of ${\cal G}_n$ is
$O(\log^2 n)$, we can verify that
$$\sum_n\sum_{(k,\ell)\in {\cal G}_n}
\mathbf{P}(\overline{C}_{n})<\infty.$$

By Borel-Cantelli lemma, with probability 1,
$\cap_{(k,\ell)\in{\cal G}_n} C_{n}(k,\ell)$
occurs for all but finitely many $n$. This completes the proof of
the statements (ii) of both Theorems 3.3 and 3.4. $\Box$

\section{Proof of Theorems 3.5 and 3.6}

\renewcommand{\thesection}{\arabic{section}} \setcounter{equation}{0}
\setcounter{theorem}{0} \setcounter{lemma}{0}

We start with the proof of Theorem 3.6 which is similar to the
proof of Theorem I. So we do not give all the details. Recall the
notations of the theorem,  Proposition 4.1  and Fact 6.

For $m_z$, given in (\ref{mz}), we have
\begg
m_z=\be(\xi(z,T)|\, \, T<\infty).
\endd

\begin{lemma} For $\log(1-2q(1-2q))<v<\log (1+2q(1-2q))$ we have
\begg \varphi(v)= \exp( m_z(v+O(v^2)),\qquad v\to 0, \label{becs0}
\endd
where  O  is uniform  in $z,$ and \begg
\psi(v)<\frac{1+|e^v-1|}{1-\frac{|e^v-1|}{1-2q}}. \label{becs1}
\endd
\end{lemma}
{\bf Proof:} The proof of this lemma is based on Proposition 4.1
and goes along the same lines as Lemma 2.3 in \cite{CsFR2}.
$\Box$

Let $k_n:=(1-\delta_n)\lambda_0\log n\sim \lambda_0\log n,$
$r_n:=c\log\log n,$ and $I(r):=[-r,r].$  Furthermore, let
$n_\ell=[e^\ell]$, $\xi(z)=\xi(z,\infty)$ and define the events
$$
A_j=\left\{\xi(S_j)\geq k_{n_\ell},\, \max_{x\in
I(r_{n_{\ell+1}})}\left(\frac{\xi(S_j+x)}
{m_xk_{n_\ell}}-1\right)\geq\varepsilon\right\}.
$$
Then
$$
\bp\left(\bigcup_{j=0}^{n_{\ell+1}}A_j\right)
\leq\sum_{j=0}^{n_{\ell+1}}\bp(A_j) \leq
\sum_{j=0}^{n_{\ell+1}}\sum_{x\in
I(r_{n_{\ell+1}})}\bp(A_j^{(x)}),
$$
where
$$
A_j^{(x)}=\left\{\xi(S_j)\geq k_{n_\ell},\, \xi(S_j+x)\geq
(1+\varepsilon)m_x k_{n_\ell}\right\}.
$$

Consider the random walk obtained by reversing the original walk
at $S_j$, i.e. let $S_i':=S_{j-i}-S_j$, $i=0,1,\ldots,j$ and
extend it to infinite time, and also the forward random walk
$S_i'':=S_{j+i}-S_j$, $i=0,1,2,\ldots$ Then $\{S_0',S_1',\ldots\}$
and $\{S_0'',S_1'',\ldots\}$ are independent random walks and so
are their respective local times $\xi'$ and $\xi''$.
Moreover,

$$
\xi(S_j)=\xi''(0)+\xi(S_j,j)\leq \xi''(0)+\xi'(0)+1,
$$
$$
\xi(S_j+x)=\xi''(x)+\xi(S_j+x,j)\leq \xi''(x)+\xi'(x).
$$
Here $\xi'$ and $\xi''$ are independent and  $\xi'$ has the same
distribution as  $\xi^* $  (see Fact 7) and $\xi''$ has the same
distribution as $\xi$.

Hence
\begin{eqnarray*}
\bp(A_j^{(x)}) &\leq&\bp(\xi''(0)+\xi'(0)\geq k_{n_\ell}-1,\,
\xi''(x)+\xi'(x)\geq
(1+\varepsilon)m_x k_{n_\ell})\\
&=&\sum \bp(\xi''(0)=k_1,\, \xi'(0)=k_2,\, \xi''(x)+\xi'(x)\geq
(1+\varepsilon) m_x k_{n_\ell}),
\end{eqnarray*}
where the summation goes for $\{(k_1,k_2):\, k_1+k_2\geq
k_{n_\ell}-1\}$. Using exponential Markov inequality, Proposition 4.1,
Fact 7, the independence of $\xi''$ and $\xi'$ and elementary calculus,
we get
\begin{eqnarray}
\bp(A_j^{(x)}) &\leq& \sum \be
\left(e^{v(\xi''(x)+\xi'(x))},\xi''(0)=k_1, \xi'(0)=k_2\right)
e^{-v(1+\varepsilon)m_x k_{n_\ell}} \nonumber \\
&=&\sum (\varphi(v))^{k_1+k_2}(1-2q)^2(2q)^{k_1+k_2}\psi(v)
e^{-v(1+\varepsilon)m_xk_{n_\ell}} \nonumber \\
&=&(1-2q)^2\psi(v)e^{-v(1+\varepsilon)m_xk_{n_\ell}}\sum
(2q\varphi(v))^{k_1+k_2}\nonumber \\
&=&(1-2q)^2\psi(v)e^{-v(1+\varepsilon)m_xk_{n_\ell}}
(2q\varphi(v))^{k_{n_\ell}} \nonumber \\
&\times&\left(\frac{k_{n_\ell}}{2q\varphi(v)(1-2q\varphi(v))}
+\frac1{(1-2q\varphi(v))^2}\right). \label{nagya}
\end{eqnarray}

Observe that even though the moment generating functions in
Proposition 4.1 and  Fact 7 are slightly different for positive
and negative values of $x$, in (\ref{nagya}) we get the same
expression while working with $\xi'+\xi''$.

By (\ref{becs0}), we obtain for all $j\geq 0$
\begin{eqnarray*}
\bp(A_j^{(x)})&\leq& (1-2q)^2\psi(v)
\left(\frac{k_{n_\ell}}{2q\varphi(v)(1-2q\varphi(v))}
+\frac1{(1-2q\varphi(v))^2}\right) \\
&\times& e^{-m_xvk_{n_\ell}(\varepsilon+O(v))}(2q)^{k_{n_\ell}}.
\end{eqnarray*}
Choose $v_0>0$ small enough such that
$$
\varepsilon+O(v_0)>0,\quad e^{v_0}<1+2q(1-2q),\quad
\frac{1}{2}<\varphi(v_0)<\frac1{2q}.
$$
Using $x\in I(r_{n_{\ell+1}})$,  we get
$$
m_xk_{n_\ell}=\frac{h^{|x|}}{2q}(1-\delta_{n_\ell})\lambda_0\log
n_\ell\geq \frac{h^{r_{n_{\ell+1}}}}{2q}
(1-\delta_{n_\ell})\lambda_0\log n_\ell.
$$
In the sequel $C_1,C_2,\ldots$ denote positive
constants whose values are unimportant in our proofs.

By the above assumptions
\begin{eqnarray*}
\bp(A_j^{(x)})&\leq& C_2k_{n_\ell}
e^{-m_xv_0 k_{n_\ell}(\varepsilon+O(v_0))}(2q)^{k_{n_\ell}} \\
&\leq& C_2k_{n_\ell} \exp \left(-(1-\delta_{n_\ell})\log n_\ell
\left(C_3h^{r_{n_{\ell+1}}} +1\right)\right).
\end{eqnarray*}
Hence
\begin{eqnarray*}
\sum_{j=0}^{n_{\ell+1}}\sum_{x\in I(r_{n_{\ell+1}})}\bp(A_j^{(x)})
&\leq& C_4n_{\ell+1}r_{n_{\ell+1}} k_{n_\ell}
\exp\left(-(1-\delta_{n_\ell})\log n_\ell
\left(C_3 h^{r_{n_{\ell+1}}}+1\right)\right) \\
&\leq& C_5\frac{n_{\ell+1}}{n_\ell}k_{n_\ell}r_{n_{\ell+1}}
\exp\left(-C_6 h^{r_{n_{\ell+1}}}\log n_\ell+
\delta_{n_\ell}\log n_\ell\right) \\
&=& C_5\frac{n_{\ell+1}}{n_\ell}k_{n_\ell}r_{n_{\ell+1}}
\exp\left(- h^{r_{n_{\ell}}}\log n_\ell
\left(C_6h^{r_{n_{\ell+1}}-r_{n_\ell}}
-\frac{\delta_{n_\ell}}{h^{r_{n_\ell}}}\right)\right) \\
&\leq& C_5\frac{n_{\ell+1}}{n_\ell}k_{n_\ell}r_{n_{\ell+1}}
\exp\left(-C_7 h^{r_{n_{\ell}}}\log n_\ell\right) \\&\leq&
C_8(\log n_\ell)\log\log
n_{\ell}\exp(-C_7(\log{n_{\ell})^{1-\alpha c}}),
\end{eqnarray*}
where in the last two lines we used the conditions of the Theorem.
Consequently,
$$
\bp(\bigcup_{j=0}^{n_{\ell+1}} A_j)\leq
\sum_{j=0}^{n_{\ell+1}}\sum_{x\in I(r_{n_{\ell+1}})}\bp(A_j^{(x)})
\leq C_8 \ell \log \ell \exp(-C_7\ell^{1-\alpha c})
$$
for large enough $\ell$, which is summable in $\ell$  when $\alpha
c<1$. By Borel-Cantelli lemma for large $\ell$ if $\xi(S_j)\geq
k_{n_\ell}$, then $\xi(S_j+x)\leq (1+\varepsilon)m_x k_{n_\ell}$
for all $x\in I(r_{n_{\ell+1}})$.

Let now $n_\ell\leq n<n_{\ell+1}$ and $x\in I(r_{n_{\ell+1}})$.
$\xi(S_j)\geq k_n, j\leq n$ implies $\xi(S_j)\geq k_{n_\ell}$,
i.e.
\begin{equation}
\xi(S_j+x)\leq (1+\varepsilon)m_xk_{n_\ell}\leq
(1+\varepsilon)m_xk_n. \label{xisx}
\end{equation}

The lower bound is similar, with slight modifications. However we
do not present it. The interested reader should look at the
corresponding proof of Theorem 1.2 in \cite{CsFR2}. ${\Box}$

\bigskip\noindent
The proof of Theorem 3.5  again goes similarly to the proof of Theorem
H. As the main ingredient is the following lemma, which is somewhat
different from the $d$-dimensional situation,  we give a complete
proof.
\begin{lemma}
Let  $0<\alpha<1$, $j\leq n-n^\alpha$, $|x|\leq c\log n$ with any $c>0$.
Then with probability 1 there exists an $n_0(\omega)$ such that for
$n\geq n_0$ we have
$$
\xi(S_j+x,n)=\xi(S_j+x,\infty).
$$
\end{lemma}

\noindent {\bf Proof.} Let
$$
A_n=\bigcup_{j\leq n-n^\alpha} \quad\bigcup_{\ell\geq n}
\quad\bigcup_{|x|\leq c\log n} \{S_\ell-S_j=x\}.
$$
By  our Lemma 4.1
\begin{equation}
\bp(S_\ell-S_j=x)=\bp(S_{\ell-j}=x)\leq C_1\exp(-C_2(\ell-j)+C_3 x).
\end{equation}

Consequently,
$$
\bp(A_n)\leq C_1\sum_{j\leq n-n^\alpha}\sum_{\ell\geq n}\sum_{|x|\leq
c\log n}\exp(-C_2\ell+C_2 j+C_3x)
$$
$$
\leq C\exp(-C_2n+C_2(n-n^\alpha)+C_3c\log n)=C_4
n^{C_5}\exp(-C_2n^\alpha).
$$

Since this is summable, we have the lemma. $\Box$

To prove Theorem 3.5, observe that it suffices to consider points
visited before time $n-n^\alpha$, ($0<\alpha<1$), since in the
time interval $(n-n^\alpha,n)$ the maximal local time is less than
$\alpha(1+\varepsilon)\lambda_0\log n$, hence this point cannot be
in ${\cal A}_n$. Consequently, Theorem 3.5 follows from Theorem
3.6 and Lemma 7.2.  $\Box$

\section{Concluding remarks}

\renewcommand{\thesection}{\arabic{section}} \setcounter{equation}{0}
\setcounter{theorem}{0} \setcounter{lemma}{0}

First observe that the following points are on the curve
$g(x,y)=1.$

\begg \left(0,\frac{1}{\log(1/q)} \right),\quad \left(
-\frac{1}{\log(2pq)},-\frac{1}{\log(2pq)} \right),\quad
\left(\frac{2\kappa_0 p}{2p+1}, \kappa_0 \right),\quad
\left(\lambda_0, \frac{\lambda_0}{p} \right).
\endd

Consequently, there are points $x_n$ such that
$$\xi(x_n,n)=1\quad {\rm and}\quad \Xi(x_n,n)\sim \frac{\log
n}{\log(1/q)}$$
which in fact means that
$$\Xi(x_n,n)=\xi(x_n+1,n ) \sim\frac{\log n}{\log(1/q)}.$$

On the other hand, if  for a point $x_n$,

$$\Xi(x_n,n)>(1+\epsilon)\frac{\log n}{\log(1/q)},$$
then we have $\xi(x_n,n)>c\log n$ with some $c>0.$

If $\xi(x_n,n)\sim \lambda_0\log n$ then for the unit sphere centered
at $x_n$, that is to say for its two neighbors we have
$$\Xi(x_n,n)\sim\frac{\lambda_0}{p}\log n.$$
Since $m_{-1}=m_1,$ it follows from Corollary 3.1 that the
two neighbors of the nearly maximally visited points have
asymptotically equal local time.

On the other hand, if the occupation time is asymptotically  maximal,
$$\Xi(x_n,n)\sim \kappa_0 \log n \quad {\rm then} \quad  \xi(x_n,n)\sim
\frac{2p}{2p+1}\kappa_0\log n.$$

With some extra calculation one can find the maximal weight of
the unit ball:

$$
w(x_n,n):=\xi(x_n,n)+\Xi(x_n,n),\quad w(n):= \max_{x_n\in
\mathcal{Z}} (\xi(x_n,n)+\Xi(x_n,n)).
$$

\noindent We get
$$
\lim_{n\to \infty}\frac{w(n)}{\log n}= \frac{2\,\beta}
{\log\left(\left(\frac{2q}{\beta+1}\right)^{\beta+1}
\left(p(\beta-1)\right)^{\beta-1}\right)}=:\frac{2\, \beta}{C},
\qquad\hspace{1cm} \mathrm{%
a.s.}  \label{last}
$$
where $\beta$ is the constant defined in (\ref{d16}). In this case
we have

$$
\xi(x_n,n)\sim\frac{\beta-1}{C}\log n,\quad  {\rm and}  \quad
\Xi(x_n,n)\sim \frac{\beta+1}{C} \log n.
$$

\noindent As a final conclusion if any  of the three quantities of
$\xi(x_n, n), \,\,\Xi(x_n, n), \,\,w(x_n, n),$ is asymptotically
maximal, it uniquely determines the asymptotic values of the other
two quantities, an interesting phenomenon which we proved for
$d>4$ in the symmetric walk  case.

\end{document}